\documentclass[10pt,reqno]{amsart}
\usepackage{graphicx}
\usepackage{indentfirst,csquotes}

\usepackage{amssymb,amsthm,amsmath}
\usepackage{xcolor,paralist,hyperref,titlesec,graphicx,fancyhdr,etoolbox,enumerate,textcomp,listings,verbatim,subfig,rotating,mathrsfs}
\newtheorem{theorem}{Theorem}[]

\newtheorem{lemma}[theorem]{Lemma}
\newtheorem{proposition}[theorem]{Proposition}

\titleformat{\section}[display]{\normalfont\huge\bfseries\centering}{\centering\chaptertitlename\thechapter}{10pt}{\Large}
\titlespacing*{\section}{0pt}{0ex}{0ex}

\hypersetup{ colorlinks=true, linkcolor=black, filecolor=black, urlcolor=black }
\def\proof{\noindent {\it Proof. $\, $}}
\def\endproof{\hfill $\Box$ \vskip 5 pt }

\definecolor{codegreen}{rgb}{0,0.6,0}
\definecolor{codegray}{rgb}{0.5,0.5,0.5}
\definecolor{codepurple}{rgb}{0.58,0,0.82}
\definecolor{backcolour}{rgb}{0.95,0.95,0.92}

\lstdefinestyle{mystyle}{
    stringstyle=\color{codepurple},
    frame=tb,
    basicstyle=\footnotesize,         
    breaklines=true,                 
    captionpos=b,                    
    keepspaces=true,                 
    showspaces=false,                
    showstringspaces=false,
    showtabs=false,                  
    tabsize=2
}

\usepackage{lipsum}

\begin{document}
\title{On a Mixed Arithmetic-Mean, Geometric-Mean, Harmonic-Mean Inequality} 
\author[K. Nam]{Kyumin Nam}
\date{\today}
\address{Incheon Posco Academy, 
50 Convensia-daero 42beon-gil, Songdo-dong, Yeonsu-gu, Incheon, South Korea}
\email{godbros.miniprime@gmail.com}
\maketitle

\let\thefootnote\relax

\begin{abstract}
A mixed arithmetic-mean, geometric-mean inequality was conjectured by F. Holland and proved by K. Kedlaya. In this note, we prove a mixed arithmetic-mean, harmonic-mean inequality and a mixed geometric-mean, harmonic-mean, and a more extended inequality: a mixed arithmetic-mean, geometric-mean, harmonic-mean inequality.
\end{abstract} 

\bigskip \noindent The following inequality, a mixed arithmetic-mean, geometric-mean inequality, was conjectured by F. Holland in \cite{1} and proved by K. Kedlaya in \cite{2}: \begin{align} \label{eq:(1)}
    \left( x_1 \cdot \frac{x_1 + x_2}{2} \cdots \frac{x_1 + x_2 + \cdots x_n}{n} \right)^{1/n} \geq \frac{1}{n} \left(x_1 + \sqrt{x_1 x_2} + \cdots + \sqrt[n]{x_1 x_2 \cdots x_n} \right)
\end{align} where $x_1,x_2,\cdots,x_n$ are positive real numbers, with equality if and only if $x_1 = x_2 = \cdots = x_n$. 

First, without proof, we will state the following lemma from \cite{2}.

\noindent \begin{lemma} The vectors ${\it {\bf a}}(i, j) = (a_1(i,j), a_2(i,j), \cdots, a_n(i,j))$ given by \begin{align*}
    a_k(i,j) &= {n-i \choose j-k} {i-1 \choose k-1} \bigg/ {n-1 \choose j-1} \\
    &= \frac{(n-i)! (n-j)! (i-1)! (j-1)!}{(n-1)! (k-1)! (n-i-j+k)! (i-k)! (j-k)!} \quad (i,j=1,2,\cdots,n)
\end{align*} satisfy \begin{enumerate}[i.]
    \item $a_k(i,j) \geq 0$ for all $i,j,k$,
    \item $a_k(i,j) = 0$ for $k > \min(i,j)$,
    \item $a_k(i,j) = a_k(j,i)$ for all $i,j,k$,
    \item $\sum_{k=1}^{n} a_k(i,j) = 1$ for all $i,j$,
    \item $\sum_{i=1}^{n} a_k(i,j) = n/j$ for $k \leq j$, and $\sum_{i=1}^{n} a_k(i,j) = 0$ for $k > j$.
\end{enumerate}
\end{lemma}

Then the proof of (1) in \cite{2} is followed as:

\noindent \begin{proposition} Let $x_1,x_2,\cdots,x_n$ be positive real numbers. Then we have \begin{align*}
        \left( \prod_{j=1}^{n} \frac{x_1 +x_2 + \cdots + x_j}{j} \right)^{1/n} \geq \frac{1}{n} \sum_{i=1}^{n} \sqrt[i]{x_1 x_2 \cdots x_i},
\end{align*} with equality if and only if $x_1 = x_2 = \cdots = x_n$.
\end{proposition}

\proof{} We use the following notation for weighted arithmetic mean and geometric means: \begin{align*}
    \mathscr{M}({\bf x}, {\bf a}) = \sum_{k=1}^{n} a_k x_k \qquad \mathrm{and} \qquad \mathscr{G}({\bf x}, {\bf a}) = \prod_{k=1}^{n} {x_k}^{a_k}
\end{align*} where ${\bf a} = (a_1,a_2,\cdots,a_n)$ is an $n$-tuple of nonnegative real numbers such that $\sum_{k=1}^{n} a_k = 1$. By the arithmetic-mean, geometric-mean inequality \begin{align}
    \mathscr{M}({\bf x}, {\bf a}) \geq \mathscr{G}({\bf x}, {\bf a})
\end{align} with equality if and only if $x_1 = x_2 = \cdots = x_n$. Let $\mathscr{M}(i, j)$ and $\mathscr{G}(i, j)$ be the means obtained by setting ${\bf a} = {\bf a}(i,j)$ in $\mathscr{M}({\bf x}, {\bf a})$ and $\mathscr{G}({\bf x}, {\bf a})$. By Lemma 1 and (2), we have \begin{align*}
    \frac{x_1 + x_2 + \cdots + x_j}{j} = \frac{1}{n} \sum_{k=1}^{n} x_k \sum_{i=1}^{n} a_k(i,j) = \frac{1}{n} \sum_{i=1}^{n} \mathscr{M}(i,j) \geq \frac{1}{n} \sum_{i=1}^{n} \mathscr{G}(i,j).
\end{align*} Taking the geometric mean of both sides over $j$, we get \begin{align}
    \left( \prod_{j=1}^{n} \frac{x_1 +x_2 + \cdots + x_j}{j} \right)^{1/n} \geq \frac{1}{n} \prod_{j=1}^{n} \left( \sum_{i=1}^{n} \mathscr{G}(i,j) \right)^{1/n}.
\end{align} By Hölder's inequality (\cite{3}, p.\,21), \begin{align}
    \frac{1}{n} \prod_{j=1}^{n} \left( \sum_{i=1}^{n} \mathscr{G}(i,j) \right)^{1/n} \geq \frac{1}{n} \sum_{i=1}^{n} \prod_{j=1}^{n} \mathscr{G}(i,j)^{1/n}.
\end{align} Equality holds only if every two ${\bf g}_1,{\bf g}_2,\cdots,{\bf g}_n$ are proportional where \begin{align*}
    {\bf g}_i = (\mathscr{G}(i,1), \mathscr{G}(i,2), \cdots, \mathscr{G}(i,n)) \quad (i=1,2,\cdots,n).
\end{align*} Since $\mathscr{G}(i,1) = x_1$ for all $i$, this would imply that ${\bf g}_1 = {\bf g}_2 = \cdots = {\bf g}_n$. And as $\mathscr{G}(i,n) = x_i$ for all $i$, this would also imply that $x_1 = x_2 = \cdots = x_n$. Here, by Lemma 1, we have \begin{align}
    \prod_{j=1}^{n} \mathscr{G}(i,j)^{1/n} = \prod_{k=1}^{n} \prod_{j=1}^{n} {x_k}^{a_k(i,j)/n} = \prod_{k=1}^{i} {x_k}^{1/i} = \sqrt[i]{x_1 x_2 \cdots x_i}.
\end{align} Combining (3), (4), and (5) completes the proof. \endproof{}

From a mixed arithmetic-mean, geometric-mean inequality, a mixed arithmetic-mean, harmonic-mean inequality can be conjectured as \begin{align*}
        \left[ \frac{1}{n} \sum_{j=1}^{n} \left(\frac{x_1 +x_2 + \cdots + x_j}{j}\right)^{-1} \right]^{-1} \geq \frac{1}{n} \sum_{i=1}^{n} \frac{i}{\frac{1}{x_1} + \frac{1}{x_2} + \cdots + \frac{1}{x_i}},
\end{align*} with equality if and only if $x_1 = x_2 = \cdots = x_n$. We can prove this inequality using a similar approach above.

\noindent \begin{proposition} Let $x_1,x_2,\cdots,x_n$ be positive real numbers. Then we have \begin{align*}
        \left[ \frac{1}{n} \sum_{j=1}^{n} \left(\frac{x_1 +x_2 + \cdots + x_j}{j}\right)^{-1} \right]^{-1} \geq \frac{1}{n} \sum_{i=1}^{n} \frac{i}{\frac{1}{x_1} + \frac{1}{x_2} + \cdots + \frac{1}{x_i}},
\end{align*} with equality if and only if $x_1 = x_2 = \cdots = x_n$.
\end{proposition}

\proof{} We use the following notation for weighted arithmetic mean and harmonic means: \begin{align*}
    \mathscr{M}({\bf x}, {\bf a}) = \sum_{k=1}^{n} a_k x_k \qquad \mathrm{and} \qquad \mathscr{H}({\bf x}, {\bf a}) = \left( \sum_{k=1}^{n} a_k {x_k}^{-1} \right)^{-1}
\end{align*} where ${\bf a} = (a_1,a_2,\cdots,a_n)$ is an $n$-tuple of nonnegative real numbers such that $\sum_{k=1}^{n} a_k = 1$. By the arithmetic-mean, harmonic-mean inequality \begin{align}
    \mathscr{M}({\bf x}, {\bf a}) \geq \mathscr{H}({\bf x}, {\bf a})
\end{align} with equality if and only if $x_1 = x_2 = \cdots = x_n$. Let $\mathscr{M}(i, j)$ and $\mathscr{H}(i, j)$ be the means obtained by setting ${\bf a} = {\bf a}(i,j)$ in $\mathscr{M}({\bf x}, {\bf a})$ and $\mathscr{H}({\bf x}, {\bf a})$. By Lemma 1 and (6), we obtain \begin{align*}
    \frac{x_1 + x_2 + \cdots + x_j}{j} = \frac{1}{n} \sum_{k=1}^{n} x_k \sum_{i=1}^{n} a_k(i,j) = \frac{1}{n} \sum_{i=1}^{n} \mathscr{M}(i,j) \geq \frac{1}{n} \sum_{i=1}^{n} \mathscr{H}(i,j).
\end{align*} Taking the harmonic mean of both sides over $j$, we get \begin{align}
    \left[ \frac{1}{n} \sum_{j=1}^{n} \left(\frac{x_1 +x_2 + \cdots + x_j}{j}\right)^{-1} \right]^{-1} 
    &\geq \left[ \frac{1}{n} \sum_{j=1}^{n} \left(\frac{1}{n} \sum_{i=1}^{n} \mathscr{H}(i,j) \right)^{-1} \right]^{-1} \\
    &= \left[\sum_{j=1}^{n} \left(\sum_{i=1}^{n} \mathscr{H}(i,j) \right)^{-1} \right]^{-1}.
\end{align} By Minkowski's inequality (\cite{3}, p.\,30), \begin{align}
    \left[\sum_{j=1}^{n} \left(\sum_{i=1}^{n} \mathscr{H}(i,j) \right)^{-1} \right]^{-1}
    &\geq \sum_{i=1}^{n} \left( \sum_{j=1}^{n} \mathscr{H}(i,j)^{-1} \right)^{-1} \\
    &=\frac{1}{n} \sum_{i=1}^{n} \left( \frac{1}{n} \sum_{j=1}^{n} \mathscr{H}(i,j)^{-1} \right)^{-1}.
\end{align} Equality holds only if ${\bf h}_1,{\bf h}_2,\cdots,{\bf h}_n$ are proportional where \begin{align*}
    {\bf h}_i = (\mathscr{H}(i,1), \mathscr{H}(i,2), \cdots, \mathscr{H}(i,n)) \quad (i=1,2,\cdots,n).
\end{align*} Since $\mathscr{H}(i,1) = x_1$ for all $i$, this would imply that ${\bf h}_1 = {\bf h}_2 = \cdots = {\bf h}_n$. And as $\mathscr{H}(i,n) = x_i$ for all $i$, this would also imply that $x_1 = x_2 = \cdots = x_n$. Here, by Lemma 1, we have \begin{align}
    \left(\frac{1}{n} \sum_{j=1}^{n} \mathscr{H}(i,j)^{-1}\right)^{-1} 
    &= \left(\frac{1}{n} \sum_{k=1}^{n} {x_k}^{-1} \sum_{j=1}^{n} a_k(i,j)\right)^{-1} \\
    &= \left(\frac{1}{i} \sum_{k=1}^{i} {x_k}^{-1}\right)^{-1} 
    = \frac{i}{\frac{1}{x_1} + \frac{1}{x_2} + \cdots + \frac{1}{x_i}}
\end{align} Combining (7), (8), (9), (10), (11), and (12) completes the proof. \endproof{}

Similarly, a mixed geometric-mean, harmonic-mean inequality also can be conjectured as \begin{align*}
    \left[ \frac{1}{n} \sum_{j=1}^{n} \left(\sqrt[j]{x_1 x_2 \cdots x_j}\right)^{-1} \right]^{-1} \geq \left(\prod_{i=1}^{n} \frac{i}{\frac{1}{x_1} + \frac{1}{x_2} + \cdots + \frac{1}{x_i}}\right)^{1/n},
\end{align*} with equality if and only if $x_1 = x_2 = \cdots = x_n$. We can also prove this inequality using a similar approach above.

\noindent \begin{proposition} Let $x_1,x_2,\cdots,x_n$ be positive real numbers. Then we have \begin{align*}
    \left[ \frac{1}{n} \sum_{j=1}^{n} \left(\sqrt[j]{x_1 x_2 \cdots x_j}\right)^{-1} \right]^{-1} \geq \left(\prod_{i=1}^{n} \frac{i}{\frac{1}{x_1} + \frac{1}{x_2} + \cdots + \frac{1}{x_i}}\right)^{1/n},
\end{align*} with equality if and only if $x_1 = x_2 = \cdots = x_n$.
\end{proposition}

\proof{} We use the following notation for weighted geometric mean and harmonic means: \begin{align*}
    \mathscr{G}({\bf x}, {\bf a}) = \prod_{k=1}^{n} {x_k}^{a_k} \qquad \mathrm{and} \qquad \mathscr{H}({\bf x}, {\bf a}) = \left( \sum_{k=1}^{n} a_k {x_k}^{-1} \right)^{-1}
\end{align*} where ${\bf a} = (a_1,a_2,\cdots,a_n)$ is an $n$-tuple of nonnegative real numbers such that $\sum_{k=1}^{n} a_k = 1$. By the geometric-mean, harmonic-mean inequality \begin{align}
    \mathscr{G}({\bf x}, {\bf a}) \geq \mathscr{H}({\bf x}, {\bf a})
\end{align} with equality if and only if $x_1 = x_2 = \cdots = x_n$. Let $\mathscr{G}(i, j)$ and $\mathscr{H}(i, j)$ be the means obtained by setting ${\bf a} = {\bf a}(i,j)$ in $\mathscr{G}({\bf x}, {\bf a})$ and $\mathscr{H}({\bf x}, {\bf a})$. By Lemma 1 and (13), we obtain \begin{align*}
    \sqrt[j]{x_1 x_2 \cdots x_j} = \left( \prod_{k=1}^{n} \prod_{i=1}^{n} {x_k}^{a_k (i,j)}\right)^{1/n}= \left(\prod_{i=1}^{n} \mathscr{G}(i, j) \right)^{1/n} \geq \left(\prod_{i=1}^{n} \mathscr{H}(i, j) \right)^{1/n}.
\end{align*} Taking the harmonic mean of both sides over $j$, we get \begin{align}
    \left[ \frac{1}{n} \sum_{j=1}^{n} \left(\sqrt[j]{x_1 x_2 \cdots x_j}\right)^{-1} \right]^{-1} &\geq \left[ \frac{1}{n} \sum_{j=1}^{n} \left\{ \left(\prod_{i=1}^{n} \mathscr{H}(i, j) \right)^{1/n} \right\}^{-1} \right]^{-1} \\
    &= \left[ \frac{1}{n} \sum_{j=1}^{n}  \left(\prod_{i=1}^{n} \mathscr{H}(i, j)^{-1} \right)^{1/n} \right]^{-1}.
\end{align} By Hölder's inequality (\cite{3}, p.\,21), \begin{align*}
    \sum_{j=1}^{n}  \left(\prod_{i=1}^{n} \mathscr{H}(i, j)^{-1} \right)^{1/n} \leq \prod_{i=1}^{n} \left( \sum_{j=1}^{n} \mathscr{H}(i,j)^{-1} \right)^{1/n},
\end{align*} therefore \begin{align}
    \left[ \frac{1}{n} \sum_{j=1}^{n}  \left(\prod_{i=1}^{n} \mathscr{H}(i, j)^{-1} \right)^{1/n} \right]^{-1} \geq \left[ \prod_{i=1}^{n} \left(\frac{1}{n} \sum_{j=1}^{n} \mathscr{H}(i, j)^{-1} \right)^{-1} \right]^{1/n}.
\end{align} Equality holds only if every two ${\bf h}_1,{\bf h}_2,\cdots,{\bf h}_n$ are proportional where \begin{align*}
    {\bf h}_i = \left(\mathscr{H}(i,1)^{-1}, \mathscr{H}(i,2)^{-1}, \cdots, \mathscr{H}(i,n)^{-1}\right) \quad (i=1,2,\cdots,n).
\end{align*} Since $\mathscr{H}(i,1)^{-1} = {x_1}^{-1}$ for all $i$, this would imply that ${\bf h}_1 = {\bf h}_2 = \cdots = {\bf h}_n$. And as $\mathscr{H}(i,n)^{-1} = {x_i}^{-1}$ for all $i$, this would also imply that $x_1 = x_2 = \cdots = x_n$. Here, by Lemma 1, we have \begin{align}
    \left(\frac{1}{n} \sum_{j=1}^{n} \mathscr{H}(i,j)^{-1}\right)^{-1} 
    &= \left(\frac{1}{n} \sum_{k=1}^{n} {x_k}^{-1} \sum_{j=1}^{n} a_k(i,j)\right)^{-1} \\
    &= \left(\frac{1}{i} \sum_{k=1}^{i} {x_k}^{-1}\right)^{-1} 
    = \frac{i}{\frac{1}{x_1} + \frac{1}{x_2} + \cdots + \frac{1}{x_i}}
\end{align} Combining (14), (15), (16), (17), and (18) completes the proof. \endproof{}

Combining these inequalities, we can get a mixed arithmetic-mean, geometric-mean, harmonic-mean inequality. Thus, we obtain the following theorem.

\noindent \begin{theorem} Let $x_1, x_2, \cdots, x_n$ be positive real numbers. Then we have \begin{align}
    &\left[ \frac{1}{n} \sum_{i=1}^{n} \left\{ \left( \prod_{j=1}^{i} \frac{x_1 + x_2 + \cdots + x_j}{j} \right)^{1/i} \right\}^{-1} \right]^{-1} \\
    &\geq \left[ \frac{1}{n} \sum_{i=1}^{n} \left( \frac{1}{i} \sum_{j=1}^{i} \sqrt[j]{x_1 x_2 \cdots x_j} \right)^{-1} \right]^{-1} \geq \frac{1}{n} \sum_{i=1}^{n} \left[ \frac{1}{i} \sum_{j=1}^{i} \left( \sqrt[j]{x_1 x_2 \cdots x_j} \right)^{-1} \right]^{-1} \\
    &\geq \frac{1}{n} \sum_{i=1}^{n} \left( \prod_{j=1}^{i} \frac{j}{\frac{1}{x_1} + \frac{1}{x_2} + \cdots + \frac{1}{x_j}} \right)^{1/i}
\end{align} with equality if and only if $x_1 = x_2 = \cdots = x_n$.
\end{theorem}

\proof{} The proof completes by applying Proposition 2, Proposition 3, and Proposition 4.  \endproof{}

$\,$


\begin{thebibliography}{99}

\bibitem{1} Holland, F. (1992). On a mixed arithmetic-mean, geometric-mean inequality, {\it Mathematics Competitions}, 5, 60-64.

\bibitem{2} Kedlaya, K. (1994). Proof of a mixed arithmetic-mean, geometric-mean inequality. {\it The American Mathematical Monthly}, 101(4), 355-357. https://doi.org/10.2307/2975630

\bibitem{3} Hardy, G., Littlewood, J. E., \& Pólya, G. (1952). {\it Inequalities} (2nd ed.). Cambridge University Press.

\end{thebibliography}
\end{document}